\documentclass[10pt]{amsart}

\usepackage{amssymb}
\usepackage{hyperref}
\usepackage{amsmath,amscd}
\usepackage[margin=1in]{geometry}
\usepackage[enableskew]{youngtab}
\usepackage{diagbox}

\newtheorem{theorem}{Theorem}[section]
\newtheorem{lemma}[theorem]{Lemma}
\newtheorem{proposition}[theorem]{Proposition}
\newtheorem{corollary}[theorem]{Corollary}

\theoremstyle{definition}

\theoremstyle{remark}
\newtheorem{remark}[theorem]{Remark}

\begin{document}

\title{Alternating Sums in Compositions} 

\author{Runde Gao}
\address{Department of Mathematics and Statistics, University of Nebraska at Kearney, NE 68849, USA}
\email{gaor@lopers.unk.edu}

\author{Jia Huang}
\address{Department of Mathematics and Statistics, University of Nebraska at Kearney, NE 68849, USA}
\email{huangj2@unk.edu}

\begin{abstract}
Motivated by the recent asymptotic formula for the moment of the alternating sum of parts in partitions by Craig and Li, we study the moment of the alternating sum of parts in compositions.
We establish closed formulas for this moment using techniques involving generating functions, matrix diagonalization, Krawtchouk polynomials, differential equations, and Fa\`a di Bruno's formula. 
Our result is connected with some known integer sequences in the OEIS and admits a combinatorial interpretation by coin flipping.
\end{abstract}



\maketitle

\section{Introduction}\label{sec:intro}

Partitions and compositions are fundamental concepts widely used in mathematics and sciences, including discrete mathematics, representation theory, symmetric function theory, statistical mechanics, genetics, etc.
See, for example, the monographs by Andrews--Eriksson~\cite{AndrewsEriksson}, Grinberg--Reiner~\cite{GrinbergReiner}, and Heubach--Mansour~\cite{HeubachMansour}.

Recall that a \emph{composition} is a sequence $\alpha=(\alpha_1,\ldots,\alpha_\ell)$ of positive integers with \emph{parts} $\alpha_1, \ldots, \alpha_\ell$, \emph{length} $\ell(\alpha):=\ell$, and \emph{size} of $|\alpha|:=\alpha_1+\cdots+\alpha_\ell$.
If $|\alpha|=n$ then we say $\alpha$ is a \emph{composition of $n$} and write $\alpha\models n$.
It is well known that the number of compositions of $n$ is $2^{n-1}$ since there is a bijection between these compositions and binary strings of length $n-1$.

A \emph{partition} $\lambda=(\lambda_1 \ge\ldots \ge \lambda_\ell >0)$ is a composition with weakly decreasing parts.
If $|\lambda|=n$ then we say $\lambda$ is a \emph{partition of $n$} and write $\lambda\vdash n$.
Although partitions are prevalent in many fields, no closed formula is known for $p(n)$, the number of partitions of $n$.
This is in sharp contrast to the aforementioned simple formula $2^{n-1}$ for the number of compositions of $n$.

On the other hand, many results on partitions have analogues for compositions. 
For example, Euler's partition theorem asserts that the number of partitions of $n$ into odd parts equals the number of partitions of $n$ into distinct parts.
Analogously, the number of compositions of $n$ with odd parts equals the number of compositions of $n+1$ with parts greater than one.
In fact, both of these numbers are equal to the \emph{Fibonacci number} $F_n$, which is defined recursively by $F_n:=n$ for $n=0,1$ and $F_n := F_{n-1}+F_{n-2}$ for $n\ge2$; see Cayley~\cite{Cayley}, Grimaldi~\cite{Grimaldi}, and Stanley~\cite[Exercise 1.35]{EC1}.
Moreover, generalizations of Euler's partition theorem by Glaisher and Franklin also have composition analogues by Munagi~\cite{Munagi} and the second author~\cite{CompRes}, respectively.
Another example is Legendre's theorem (also a consequence of Euler's pentagonal number theorem) on the number of even partitions of $n$ into distinct parts minus the number of odd partitions of $n$ into distinct parts, which admits a composition analogue by the second author~\cite{CompEvenOdd}.

This paper is motivated by studies on $a(\lambda):=\lambda_1-\lambda_2+\lambda_3-\cdots+(-1)^{\ell-1} \lambda_\ell$, the alternating sum of parts of a partition $\lambda=(\lambda_1,\ldots,\lambda_\ell)$.
It is known that $a(\lambda)$ counts odd parts in the conjugate partition $\lambda^t$.
For instance, $(54331)^t = 54421$ as shown by the Young diagram below, and its number of odd parts equals $a(54331) = 5-4+3-3+1=2$.
\begin{center}
\footnotesize
\yng(5,4,3,3,1)   
\qquad
\yng(5,4,4,2,1)
\end{center}
Recently, Craig and Li~\cite{CraigLi} established the following asymptotic formula for the \emph{$m$th moment} of $a(\lambda)$:
\[E_m(n) := \displaystyle\sum_{\lambda\vdash n} \frac{a(\lambda)^m}{p(n)}
\sim \frac{(6n)^{m/2}}{(2\pi)^m} \log^m \left(\frac{\sqrt{6n}}{\pi}\right).
\]

It is natural to investigate $a(\alpha):=\alpha_1-\alpha_2+\alpha_3-\cdots+(-1)^{\ell-1} \alpha_\ell$, the alternating sum of parts of a composition $\alpha=(\alpha_1,\ldots,\alpha_\ell)$, and its \emph{$m$th moment}
\[ 
E_m^c(n) := \displaystyle\sum_{\alpha\models n} \frac{a(\alpha)^m}{2^{n-1}}.
\]
Table~\ref{tab1} shows $a(\alpha)$ for all compositions $\alpha\models 4$ with mean $E^c_1(4) = 1$ and second moment $E^c_2(4)=4$ (hence variance $E^c_2(4) - E^c_1(4)^2 = 3$). 

\begin{table}[ht]
\begin{center}
\begin{tabular}{c|cccccccc}
$\alpha$ & 4 & 31 & 13 & 22 & 211 & 121 & 112 & 1111 \\
\hline
$a(\alpha)$ & 4 & 2 & -2 & 0 & 2 & 0 & 2 & 0 \\
\end{tabular}
\vskip5pt
\caption{Alternating sum of parts in compositions of $n=4$}
\label{tab1}
\end{center}
\end{table}



In an undergraduate research project mentored by the second author, Kitamura~\cite{Kitamura} computed $E_m^c(n)$ for small values of $m$ and $n$ with result suggesting that $E_m^c(n)$ is not only always a positive integer but also a polynomial in $n$:
\begin{align*}
E_1^c(n) &= 1, \\
E_2^c(n) &= n, \\ 
E_3^c(n) &= 3n-2, \\ 
E_4^c(n) &= 3n^2-2n, \\ 
E_5^c(n) &= 15n^2-30n+16, \\ 
E_6^c(n) &= 15n^3 - 30n^2 + 16n, \\ 
E_7^c(n) &= 105n^3-420n^2+588n-272,\\ 
E_8^c(n) &= 105n^4-420n^3+588n^2-272n, \\ 
E_9^c(n) &= 945 n^4 - 6300 n^3 + 16380 n^2 - 18960 n + 7936, \\
E_{10}^c(n) &= 945 n^5 - 6300 n^4 + 16380 n^3 - 18960 n^2 + 7936n, \\
E_{11}^c(n) &= 10395 n^6 - 103950 n^5 + 429660 n^4 - 893640 n^3 + 911328n^2 - 353792.
\end{align*}

Now we confirm the above polynomial expressions of $E_m^c(n)$ for small values of $m$ by establishing two closed formulas of $E_m^c(n)$ for all $m\ge0$ and $n\ge1$:

\begin{theorem}\label{thm1}
For $m\ge0$ and $n\ge1$ we have
\begin{equation}\label{eq:E1}
E^c_m(n) = \frac{1}{2^{n-1}} \sum_{i=0}^{n-1} \binom{n-1}{i} (n-2i)^m. 
\end{equation}
\end{theorem}

\begin{theorem}\label{thm2}
For $m\ge0$ and $n\ge1$ we have
\begin{equation}\label{eq:E2}
E^c_m(n) = \sum_{k=0}^m S(m+1,k+1) \sum_{ \substack{ \pi\in\Pi_k \\ B\in \pi \Rightarrow |B|>1}} \frac{(n-1)_{|\pi|}}{2^{|\pi|}} \prod_{B\in \pi} (-1)^{|B|} |B|!. 
\end{equation}
Here $\Pi_m$ is the set of all partitions of $\{1,2,\ldots,m\}$, $|\pi|$ is the number of blocks of $\pi$, $|B|$ is the size of the block $B\in \pi$, and $(n-1)_{|\pi|} := (n-1)(n-2)\cdots(n-|\pi|)$ is the falling factorial.
\end{theorem}

In particular, the mean and variance of $a(\alpha)$ for all $\alpha\models n$ are $E_1^c(n)=1$ and $E_2^c(n) - E_1^c(n)^2 = n-1$, respectively.
To prove the above theorems, we first obtain a generating function
\begin{equation}\label{eq:x^a}
\sum_{\alpha\models n} x^{a(\alpha)}
= x(x+x^{-1})^{n-1}.
\end{equation}
Then we apply the operator $\Omega:=x\frac{d}{dx}$ exactly $m$ times and derive the result in two ways, one using matrix diagonalization and the other using differential equations, both giving the first formula~\eqref{eq:E1} of $E^c_m(n)$.
For the second formula~\eqref{eq:E2}, we apply a slightly different operator $\Psi:=\frac{d}{dx} x$ and use the chain rule for higher order derivatives, which is known as Fa\`a di Bruno's formula.

We also show that the exponential generating function of $E_m^c(n)$ is
\begin{equation}\label{eq:egf}
\sum_{m=0}^\infty \frac{E_m^c(n) t^m}{m!} = e^t \cosh^{n-1}(t).
\end{equation}
Separating the even and odd powers in the above generating function, we obtain the following result.

\begin{theorem}\label{thm3}
For $m\ge1$, we have
\[
E_{2m}^c(n) = \sum_{i=1}^m (-1)^{m-i} T_{m-1,i-1} n^i
= n E_{2m-1}^c(n),
\]
where $T_{m,i}$ is the sequence~\cite[A085734]{OEIS} with exponential generating function $\sec^n(t)$.
\end{theorem}

In particular, $E_{2m}^c(n)$ is a polynomial in $n$ of degree $m$ with leading term given by the double factorial and constant term given by the \emph{tangent number}~\cite[A000182]{OEIS}.
It would be interesting to explore the connections with various interpretations of the sequence $T_{m,i}$~\cite[A085734]{OEIS}, such as the following:
\begin{itemize}
\item
As a polynomial of $n$, $E_{2m}^c(n)$ coincides with the zeta polynomial for the poset of even-sized subsets of $\{1,2,\ldots,2m\}$ ordered by inclusion (see Remark~\ref{rem:T}).
\item
The number $T_{m,i}$ enumerates down-up permutations on $\{1,2,\ldots,2m+2\}$ with $i+1$ left-to-right maxima. 
\end{itemize}
It might be worth checking some closely related sequences~\cite[A083061, A094665]{OEIS} as well.

Finally, we provide a combinatorial interpretation for $E_m^c(n)$ using a coin-flipping model. 

\begin{theorem}\label{thm4}
If we begin with $n-1$ coins all showing tails and flip at most one of them at each step, then the number of ways to have all tails right after $m$ steps is $E^c_m(n)$.
\end{theorem}

This paper is structured as follows.
In Section~\ref{sec:gen} we prove the generating function~\eqref{eq:x^a} for the alternating sum of parts in compositions and show that applying certain differential operators leads to $E^c_m(n)$.
We then derive the first formula~\eqref{eq:E1} of $E^c_m(n)$ given in Theorem~\ref{thm1} by applying the operator $\Omega:=x\frac{d}{dx}$ exactly $m$ times to Eq.~\eqref{eq:x^a} and simplify the result in two ways: by matrix diagonalization in Section~\ref{sec:matrix} and by differential equations in Section~\ref{sec:DE}.
In Section~\ref{sec:Bruno} we use the operator $\Psi:=\frac{d}{dx} x$ instead of $\Omega$ to establish the second formula~\eqref{eq:E2} of $E_m^c(n)$ appearing in Theorem~\ref{thm2} with the help of Fa\`a di Bruno's formula.
In Section~\ref{sec:egf} we obtain the exponential generating function~\eqref{eq:egf} of $E_m^c(n)$ and prove Theorem~\ref{thm3}.
In Section~\ref{sec:interpretation} we prove Theorem~\ref{thm4}, which provides a combinatorial interpretation for $E_m^c(n)$ by flipping $n-1$ coins.

\section{Generating function}\label{sec:gen}
In this section we establish the generating function~\eqref{eq:x^a} for the alternating sum of parts in compositions, from which we can derive $E_m^c(n)$.

\begin{proposition}
For $n\ge1$, we have
\[
\sum_{\alpha\models n} x^{a(\alpha)} = x \left( x+x^{-1} \right)^{n-1}.
\]
\end{proposition}

\begin{proof}[Analytic Proof]
We can write down the generating function for compositions using a variable $t$ to keep track of the parts and another variable $x$ to keep track of the alternating sum of parts.
For compositions with even lengths we have
\[
\sum_{n=0}^\infty \sum_{ \substack{\alpha\models n \\ \ell(\alpha) = 2\ell}} x^{a(\alpha)} t^n
=\sum_{\ell=0}^{\infty}\left(xt+x^2t^2+\cdots\right)^\ell \left(x^{-1}t+x^{-2}t^2+\cdots\right)^\ell.
\]
For compositions with odd lengths we have
\[
\sum_{n=0}^\infty \sum_{ \substack{\alpha\models n \\ \ell(\alpha) = 2\ell+1}} x^{a(\alpha)} t^n
=\sum_{\ell=0}^{\infty}\left(xt+x^2t^2+\cdots\right)^\ell \left(x^{-1}t+x^{-2}t^2+\cdots\right)^\ell \left(xt+x^2t^2+\cdots\right).
\]
It follows that
\begin{align*}
\sum_{n=0}^\infty \sum_{\alpha\models n} x^{a(\alpha)} t^n
&=\sum_{\ell=0}^{\infty}\left(xt+x^2t^2+\cdots\right)^\ell \left(x^{-1}t+x^{-2}t^2+\cdots\right)^\ell \left(1+xt+x^2t^2+\cdots\right) \\
&= \left(1-\frac{xt}{1-xt}\frac{x^{-1}t}{1-x^{-1}t}\right)^{-1} (1-xt)^{-1} \\
&= \left(1-xt-\frac{t^2}{1-x^{-1}t}\right)^{-1} \\
&= \frac{1-x^{-1}t}{1-xt-x^{-1}t} \\
&= 1 + \frac{xt}{1-xt-x^{-1}t} \\
&= 1 + xt \sum_{n=0}^\infty \left(xt+x^{-1}t \right)^n.
\end{align*}
Extracting the coefficient of $t^n$ gives the desired result.
\end{proof}

\begin{proof}[Combinatorial Proof]
A compositions $\alpha\models n$ corresponds to a row of $n$ boxes with the first box colored black and each other box colored either black or white.
The alternating sum of parts of $\alpha$ equals the number of black boxes minus the number of white boxes.
The result follows immediately.
\end{proof}

Now we can derive $E_m^c(n)$ by applying certain differential operators.

\begin{corollary}\label{cor:E^c}
Define $\Omega := x\frac{d}{dx}$, $\Psi := \frac{d}{dx} x$, and $u:=x+x^{-1}$.
Then
\[
E^c_m(n) = \frac{1}{2^{n-1}}\Omega^m(xu^{n-1}) \big|_{x=1}
= \frac{1}{2^{n-1}} \Psi^m(u^{n-1})|_{x=1}.
\]
\end{corollary}

\begin{proof}
Applying the operator $\Omega$ exactly $m$ times to Eq.~\eqref{eq:x^a} and substituting in $x=1$ we obtain
\[
\sum_{\alpha\models n} a(\alpha)^m = \Omega^m(xu^{n-1}) \big|_{x=1}.
\]
Divide both sides by $2^{n-1}$ gives $E^c_m(n)$.
It is similar for the operator $\Psi$.
\end{proof}

\section{Proof of Theorem~\ref{thm1} by matrix diagonalization}\label{sec:matrix}
In this section we derive the closed formula~\eqref{eq:E1} of $E^c_m(n)$ from Corollary~\ref{cor:E^c} using matrix diagonalization.

We first define $v := \Omega u = x - x^{-1}$.
Then $\Omega v = u$ and
\begin{align*} 
\Omega(xu^{n-1}) &= xu^{n-1} + (n-1) xu^{n-2}v, \\
\Omega(xu^{n-2}v) &= xu^{n-1} + xu^{n-2}v + (n-2)xu^{n-3}v^2, \\
\Omega(xu^{n-3}v^2) &= 2xu^{n-2}v + xu^{n-3}v^2 + (n-3)xu^{n-4}v^3, \\
\cdots & \qquad \cdots \\
\Omega(xuv^{n-2}) &= (n-2)xu^2 v^{n-3} + xuv^{n-2} + xv^{n-1}, \\
\Omega(xv^{n-1}) &= (n-1)xu v^{n-2} + xv^{n-1}.
\end{align*}
This is summarized in the proposition below, which is routine to prove.

\begin{proposition}\label{prop:Omega}
For $k=0,1,\ldots,n-1$, we have
\[
\Omega(xu^{n-k-1}v^k) = 
kxu^{n-k}v^{k-1} + xu^{n-k-1}v^k + (n-k-1)xu^{n-k-2}v^{k+1}.
\]
\end{proposition}

A lemma follows immediately.

\begin{lemma}\label{lem:a}
For $m\ge0$, we have
\[
\Omega^m(xu^{n-1}) = \sum_{k=0}^{n-1} a_{m,k} xu^{n-k-1}v^k
\]
where $[a_{m,0}, a_{m,1}, \ldots, a_{m,n-1}]^T =  A^m [1, 0, \ldots, 0]^T$ with
$A=[A_{kj}]_{k,j=0}^{n-1}$ defined by  
\[ 
A_{kj} := 
\begin{cases}
1, & \text{if } j=k; \\
k+1, & \text{if } j=k+1; \\
n-k, & \text{if } j=k-1; \\
0, & \text{otherwise}.
\end{cases}
\]
\end{lemma}

The matrix $A$ has eigenvalues $\lambda_i := n-2i$ for $i=0,1,\ldots,n-1$, and the eigenvector associated with $\lambda_i$ consists of the \emph{Krawtchouk polynomial} 
\[ 
P_{k,i} = P_{k,i}(n-1) := \sum_{j=0}^k (-1)^j \binom{i}{j} \binom{n-1-i}{k-j}.
\]
In fact, the following recursive relation of the Krawtchouk polynomials~\cite[Theorem~4.1]{Coleman}
\[ 
(k+1) P_{k+1,i} = (n-1-2i)P_{k,i} - (n-k)P_{k-1,i}
\]
implies that $AP = P\Lambda$, where $P := [P_{k,i}]_{k,i=0}^{n-1}$ and $\Lambda=\mathrm{diag}(\lambda_i)_{i=0}^{n-1}$.
For more properties of the Krawtchouk polynomials, see Coleman~\cite{Coleman} (note that our $P_{k,i}(n-1)$ is their $K_{k, n-1, 2}(i)$).

We can also find $Q := P^{-1}$ using the orthogonality of the Krawtchouk polynomials~\cite[Prop.~3.1]{Coleman}:
\[
\sum_{i=0}^{n-1} \binom{n-1}{i} P_{k,i} P_{j,i} 
= 2^{n-1} \binom{n-1}{j} \delta_{j,k}.
\]
Here $\delta_{i,j}$ is the Kronecker delta.
We have
\[
Q_{i,k} = \frac{1}{2^{n-1}} \binom{n-1}{i} \binom{n-1}{k}^{-1} P_{k,i}
\] 
since
\[
(PQ)_{kj} = \sum_{i=0}^{n-1} P_{ki} Q_{ij} 
= \sum_{i=0}^{n-1} \frac{1}{2^{n-1}} \binom{n-1}{i} \binom{n-1}{j}^{-1} P_{ki} P_{ji} 
= \delta_{jk}.
\]
In summary, we have the following result.

\begin{proposition}
The diagonalization of the matrix $A$ is given by $A = P \Lambda Q$.
\end{proposition}

Now we can apply $A^m = P\Lambda^m Q$ to the initial vector $[1,0,\ldots,0]^T$ to find $a_{m,k}$.

\begin{theorem}\label{thm:a}
For $m\ge0$ and $k=0,1,\ldots,n-1$ we have
\[
a_{m,k} = \frac{1}{2^{n-1}} \sum_{i=0}^{n-1} \binom{n-1}{i} (n-2i)^m \sum_{j=0}^k (-1)^j \binom{i}{j} \binom{n-1-i}{k-j}.
\]
\end{theorem}

\begin{proof}
This follows immediately from $a_{m,k} = P_{k,i} Q_{i,0} \lambda_i^m $, the $k$th component of $A^m[1,0,\ldots,0]^T$.
\end{proof}

Now we are ready to provide the first proof of Theorem~\ref{thm1}.

\begin{proof}[First proof of Theorem~\ref{thm1}]
Since $u|_{x=1}=2$ and $v|_{x=1}=0$, it follows from Corollary~\ref{cor:E^c} and Lemma~\ref{lem:a} that 
\[
E^c_m(n) =  \frac{1}{2^{n-1}}\Omega^m(xu^{n-1})|_{x=1}
= a_{m,0}.
\]
Thus taking $k=0$ in Theorem~\ref{thm:a} gives Theorem~\ref{thm1}.
\end{proof}

\begin{remark}\label{rem:d}
Computations suggests that $a_{m,k}$ involves the \emph{falling factorial} $(n-1)_k:=(n-1)(n-2)\cdots(n-k)$.
In fact, we can define $[b_{m,0}, b_{m,1}, \ldots, b_{m,n-1}]^T = B^m [1, 0, \ldots, 0]^T$ where $B = [B_{k,j}]_{k,j=0}^{n-1}$ is defined by
\[ 
B_{kj} := 
\begin{cases}
1, & \text{if } k=j \text{ or } k=j+1; \\
(k+1)(n-k-1), & \text{if } k=j-1; \\
0, & \text{otherwise}.
\end{cases}
\]
One can check that $B = D A D^{-1}$, where $D:=\mathrm{diag}((n-k-1)!)_{k=0}^{n-1}$.
It follows that $BDP = DAP = DP\Lambda$.
Thus we can apply $B^m = DP\Lambda^m QD^{-1}$ to the initial vector $[1,0,\ldots,0]^T$ and obtain
\begin{align*}
b_{m,k} &= \sum_{i=0}^m (DP)_{k,i} (QD^{-1})_{i,0} \lambda_i^m 
= \frac{a_{m,k}}{(n-1)_k} \\
&= \frac{1}{2^{n-1}(n-1)_k} \sum_{i=0}^{n-1} \binom{n-1}{i} (n-2i)^m \sum_{j=0}^k (-1)^j \binom{i}{j} \binom{n-1-i}{k-j}.
\end{align*}
Setting $k=0$ gives the same formula for $E_m^c(n) = b_{m,0}$ as in Theorem~\ref{thm1}.
\end{remark}

\section{Proof of Theorem~\ref{thm1} by differential equations}\label{sec:DE}
We can also use differential equation techniques to derive the closed formula~\eqref{eq:E1} of $E_m^c(n)$ in Theorem~\ref{thm1}. 

\begin{lemma}\label{lem:c}
Define $c_{0,0} := 1$, $c_{m,k} := 0$ whenever $k<0$ or $m<k$, and 
\[ 
c_{m,k} := (n-k) c_{m-1,k-1} + c_{m-1,k} + (k+1) c_{m-1,k+1}
\]
for $m\ge1$ and $0\le k\le m$.
Then
\[
\Omega^m(xu^{n-1}) = \sum_{k=0}^m c_{m,k} xu^{n-k-1}v^k.
\]
\end{lemma}

\begin{proof}
We prove the result by induction on $m$.
It is trivial when $m=0$.
Assume the result holds for $m\ge0$.
This together with Proposition~\ref{prop:Omega} implies
\begin{align*}
\Omega^{m+1}(xu^{n-1}) &= \sum_{k=0}^m c_{m,k} \Omega(xu^{n-k-1}v^k ) \\
&= \sum_{k=0}^m c_{m,k} (kxu^{n-k}v^{k-1} + xu^{n-k-1}v^k + (n-k-1)xu^{n-k-2}v^{k+1}) \\
&= \sum_{k=0}^{m+1} ( (n-k)c_{m,k-1} + c_{m,k} + (k+1) c_{m,k+1} )  xu^{n-k-1}v^k \\
&= \sum_{k=0}^{m+1} c_{m+1,k} xu^{n-k-1}v^k. \qedhere
\end{align*}
\end{proof}

We use differential equations to derive a closed formula for $c_{m,k}$, which agrees with the formula of $a_{m,k}$ in Theorem~\ref{thm:a}.

\begin{theorem}\label{thm:c}
For $m\ge0$ and $0\le k\le m$ we have
\[
c_{m,k} = \frac{1}{2^{n-1}} \sum_{i=0}^{n-1} \binom{n-1}{i} (n-2i)^m \sum_{j=0}^k (-1)^j \binom{i}{j} \binom{n-1-i}{k-j}.
\] 
\end{theorem}

\begin{proof}
Let $C_m(z) := \sum_{k=0}^m c(m,k) z^k$.
We have  
\begin{align*}
C_m(z) &= \sum_k (n-k)c_{m-1,k-1} z^k + \sum_k c_{m-1,k} z^k + \sum_k (k+1)c_{m-1,k+1} z^k \\
&= (n-1)zC_{m-1}(z) - z^2 C'_{m-1}(z) + C_{m-1}(z) + C'_{m-1}(z) \\
&= (1+(n-1)z)C_{m-1}(z) + (1-z^2)C'_{m-1}(z).
\end{align*}
The differential operator $L := 1 + (n-1)z + (1-z^2)\frac{d}{dz}$ has eigenvalues $\lambda_i := n-2i$ and corresponding eigenfunctions $w_i(z) := (1-z)^i (1+x)^{n-1-i}$ for $i=0,1,\ldots,n-1$.
This can be obtained by solving the differential equation $Lw = \lambda w$, which is equivalent to
\[ 
(1-z^2) \frac{dw}{dz} = (\lambda-1-(n-1)z)w. 
\]
Separating $w$ and $z$ gives
\[ 
\frac{1}{w}dw = \frac{\lambda-1-(n-1)z}{1-z^2} dz.
\]
Using partial fractions we obtain
\[ 
\int \frac{1}{w}dw = \int \left( \frac{\lambda-n}{2(1-z)} + \frac{\lambda+n-2}{2(1+z)} \right) dz.
\]
With $i:=-(\lambda-n)/2$, we have
\[
\ln w = i\ln(1-z) + (n-1-i) \ln(1+z) + c.
\] 
Exponentiate both sides gives the eigenfunctions $w_i = (1-z)^i (1+z)^{n-1-i}$ for $i=0,1,\ldots,n-1$.

Now we expand $C_0(z) = 1 = (1-z+1+z)/2$ in terms of these eigenfunctions:
\begin{align*}
1 &= \frac1{2^{n-1}} (1-z + 1+z)^{n-1} \\
&= \frac{1}{2^{n-1}} \sum_{i=0}^{n-1} \binom{n-1}{i} (1-z)^i (1+z)^{n-1-i}.
\end{align*}
Applying $L^m$ amounts to multiplying each eigenfunction by the $m$th power of the associated eigenvalue:
\[
C_m(z) = L^m(1) = \frac{1}{2^{n-1}} \sum_{i=0}^{n-1} \binom{n-1}{i} (n-2i)^m (1-z)^i (1+z)^{n-1-i}.
\]
Extracting the coefficient of $z^k$ gives the desired formula for $c_{m,k}$.
\end{proof}

Now we present another proof of Theorem~\ref{thm1}.

\begin{proof}[Second Proof of Theorem~\ref{thm1}]
We have $E_m^c(n) = \frac{1}{2^{n-1}}\Omega^m(xu^{n-1})|_{x=1} = c_{m,0}$ by Corollary~\ref{cor:E^c} and Lemma~\ref{lem:c}.
Setting $k=0$ in Theorem~\ref{thm:c} proves Theorem~\ref{thm1}.
\end{proof}


\begin{remark}
One can shows by induction on $m$ that
 \[
\Omega^m(xu^{n-1}) = \sum_{k=0}^m (n-1)_k d_{m,k} xu^{n-k-1}v^k,
\]
where $d_{0,0} := 1$, $d_{m,k} := 0$ whenever $k<0$ or $m<k$, and 
\[ d_{m,k} := d_{m-1,k-1} + d_{m-1,k} + (k+1)(n-k-1) d_{m-1,k+1}
\]
for $m\ge1$ and $0\le k\le m$.
Setting $k=0$ gives the same formula for $E_m^c(n) = d_{m,0}$ as in Theorem~\ref{thm1}.
It might be interesting to study $d_{m,k}$ as a polynomial in $n$ for $k\ge1$ (see Table~\ref{tab2} for small values of $m$ and $k$), although we cannot find its coefficients in OEIS~\cite{OEIS}.

\begin{table}[h]
\tiny
\begin{tabular}{c|ccccccccc}
\diagbox{$m$}{$k$}& $0$ & $1$ & $2$ & $3$ & $4$ & $5$ \\
\hline
$0$ & $1$ & \\ 
$1$ &  $1$ & $1$ \\ 
$2$ &  $n$ & $2$ & $1$ \\ 
$3$ &  $3n-2$ & $3n-2$ & $3$ & $1$ \\ 
$4$ &  $3n^2-2n$ & $12n-16$ & $6n-8$ & $4$ & $1$ \\ 
$5$ &  $15n^2-30n+16$ & $15n^2-30n+16$ & $30n-60$ & $10n-20$ & $5$ & $1$ \\ 
$6$ & $15n^3-30n^2+16n$ & $90n^2-300n+272$ & $45n^2-150n+136$ & $60n-160$ & $15n-40$ & $6$ & $1$ 
\end{tabular}
\vskip5pt
\caption{$d_{m,k}$ for small values of $m$ and $k$}
\label{tab2}
\end{table}
\end{remark}

\begin{remark}
One may notice that Lemma~\ref{lem:a} and Lemma~\ref{lem:c} are identical.
In fact, for any given positive integer $n$, we have $c_{m,k} = 0$ for $k\ge n$ by either the definition or the formula of $c_{m,k}$, and $a_{m,k} = c_{m,k}$ for $k=0,1,\ldots,n-1$.
Similarly, $b_{m,k}=d_{m,k}$ for $k=0,1,\ldots,n-1$ and $d_{m,k}=0$ for $k\ge n$.
\end{remark}

\section{Proof of Theorem~\ref{thm2} by Fa\`a di Bruno's formula}\label{sec:Bruno}
In this section we derive another formula for $E^c_m(n)$ using the following chain rule for higher order derivatives, which is known as Fa\`a di Bruno's formula.

\begin{theorem}[Fa\`a di Bruno's formula]
Applying $\frac{d}{dx}$ to a composite function $f\circ g$ exactly $m$ times gives the following for all $m\ge1$ (assuming the existence of all necessary derivatives of $f$ and $g$):
\[
\frac{d^m}{dx^m} (f\circ g) = \sum_{\pi\in\Pi_m} \left( f^{(|\pi|)}\circ g \right) \prod_{B\in \pi} g^{|B|}.
\] 
Here $\Pi_m$ is the set of all partitions of $\{1,2,\ldots,m\}$, $|\pi|$ is the number of blocks of $\pi$, and $|B|$ is the size of the block $B\in \pi$.
\end{theorem}

We first apply the operator $\Psi := \frac{d}{dx} x$ repeatedly to $h(x) := (x+x^{-1})^{n-1}$ and obtain
\begin{align*}
\Psi(h(x)) &= h(x) + x h'(x), \\
\Psi^2(h(x)) &= h(x) + xh'(x) + 2xh'(x) + x^2h''(x) \\
&= h(x) + 3xh'(x) + x^2h''(x), \\
\Psi^3(h(x)) &= h(x) + xh'(x) + 6x h'(x) + 3x^2 h''(x) + 3x^2 h''(x) + x^3 h'''(x) \\
&= h(x) + 7xh'(x) + 6x^2h''(x) + x^3h'''(x),
\end{align*}
and so on. 
In general, we have the following result.

\begin{proposition}\label{prop:Psi}
Let $S(m,k)$ be the Stirling number of the second kind. 
Then
\[
\Psi^m (h(x)) = \sum_{k=0}^m S(m+1,k+1) x^k h^{(k)}(x).
\]
\end{proposition}

\begin{proof}
The result is trivial when $m=0$. 
For $m\ge1$, using the inductive hypothesis we have
\begin{align*}
\Psi^m (h(x)) &= \sum_{k=0}^{m-1} S(m,k+1) \Psi \left( x^k h^{(k)}(x) \right) \\
&= \sum_{k=0}^{m-1} S(m,k+1) \left( (k+1)x^k h^{(k)}(x) + x^{k+1} h^{(k+1)}(x) \right) \\
&= \sum_{k=0}^{m} \left( S(m,k) + (k+1)S(m,k+1) \right) x^k h^{(k)}(x) \\
&= \sum_{k=0}^{m} S(m+1,k+1) x^k h^{(k)}(x). \qedhere
\end{align*}
\end{proof}

Now we provide a proof for Theorem~\ref{thm2}, which is rephrased below.

\begin{theorem}
Using the falling factorial $(n-k):=(n-1)(n-2)\cdots(n-k)$, we have 
\begin{align*}
E^c_m(n) 
&= \sum_{k=0}^m S(m+1,k+1) \sum_{ \substack{ \pi\in\Pi_k \\ B\in \pi \Rightarrow |B|>1}} \frac{(n-1)_{|\pi|}}{2^{|\pi|}} \prod_{B\in \pi} (-1)^{|B|} |B|!. \\
&
\end{align*}
\end{theorem}

\begin{proof}
By Corollary~\ref{cor:E^c}, $E^c_m(n) = \frac1{2^{n-1}}\Psi^m(h(x))|_{x=1}$.
We use Proposition~\ref{prop:Psi} to calculate $\Psi^m(h(x))$.
Since $h(x) = f(g(x))$ with $f(u) = u^{n-1}$ and $u = g(x) = x + x^{-1}$, we can apply Fa\`a di Bruno's formula to $h^{(k)}(x)$.
We have $g'(x) = 1-x^{-2}$, $g''(x) = 2x^{-3}$, $g'''(x) = -2\cdot3x^{-4}$, and in general, \
\[
g^{(k)}(x) = (-1)^k k! x^{-k-1} \quad\text{ for $k\ge2$}.
\]
Thus $g'(1)=0$ and $g^{(k)}(1) = (-1)^k k!$ for $k\ge2$.
We also have
$f^{(k)} (u) = (n-1)(n-2) \cdots (n-k) u^{n-k-1}$, which becomes $(n-1)_k 2^{n-k-1}$ when $x=1$.
The desired formula for $E_m^c(n)$ follows.
\end{proof}

\section{Exponential Generating function of $E^c_m(n)$}\label{sec:egf}
To better understand $E^c_m(n)$, we study its exponential generating function.

\begin{corollary}
The exponential generating function of $E_m^c(n)$ is
\[
\sum_{m=0}^\infty \frac{E_m^c(n) t^m}{m!} = e^t \cosh^{n-1}(t).
\]
\end{corollary}

\begin{proof}
By Theorem~\ref{thm1} we have
\begin{align*}
\sum_{m=0}^\infty \frac{E_m^c(n) t^m}{m!}
&= \frac{1}{2^{n-1}} \sum_{m=0}^\infty \sum_{i=0}^{n-1} \binom{n-1}{i} \frac{(n-2i)^m t^m}{m!} \\
&= \frac{1}{2^{n-1}} \sum_{i=0}^{n-1} \binom{n-1}{i} \sum_{m=0}^\infty \frac{(n-2i)^m t^m}{m!} \\
&= \frac{1}{2^{n-1}} \sum_{i=0}^{n-1} \binom{n-1}{i} e^{(n-2i)t} \\
&= \frac{e^{nt}}{2^{n-1}} (1+e^{-2t})^{n-1} \\
&= e^t \cosh^{n-1}(t). \qedhere
\end{align*}
\end{proof}

Now we can identify the coefficients of $E_m^c(n)$ as a polynomial in $n$ with a sequence $T_{m,k}$~\cite[A085734]{OEIS} defined by $T_{0,0} := 1$ and 
\[ 
T_{m,k} := \sum_{j=0}^{m-k} 2^j \binom{k+j}{j+1} + \binom{k+j+1}{j+1} T_{m-1,k+j-1}.
\]
For small values of $m$ and $k$, see Table~\ref{tab3}.

\begin{table}[h]
\small
\begin{tabular}{c|ccccccccc}
\diagbox{$m$}{$i$}& $0$ & $1$ & $2$ & $3$ & $4$ & $5$ \\
\hline
$0$ & $1$ & \\ 
$1$ & $2$ & $3$ \\ 
$2$ & $16$ & $30$ & $15$ \\ 
$3$ & $272$ & $588$ & $420$ & $105$ \\ 
$4$ & $7936$ & $18960$ & $16380$ & $6300$ & $945$\\
$5$ & $353792$ & $911328$ & $893640$ & $429660$ & $103950$ & $10395$
\end{tabular}
\vskip5pt
\caption{$T_{m,k}$ for small values of $m$ and $k$}
\label{tab3}
\end{table}

\begin{corollary}\label{cor:T}
For $m\ge1$, we have
\[ 
E_{2m}^c(n) = \sum_{i=1}^m (-1)^{m-i} T_{m-1,i-1} n^i
= n E_{2m-1}^c(n).
\]
\end{corollary}

\begin{proof}
We have $e^t \cosh^{n-1}(t) =\cosh^n(t) + \sinh(t)\cosh^{n-1}(t)$.
Thus the even powers of $t$ in $e^t \cosh^{n-1}(t)$ are given by $\cosh^n(t)$, that is,
\[ 
\sum_{m=0}^\infty \frac{E_{2m}^c(n) t^{2m}}{(2m)!} = \cosh^n(t).
\]
Substituting in $-n$ for $n$ and $-t^2$ for $t^2$ we have
\[ 
\sum_{m=0}^\infty \frac{E_{2m}^c(-n) (-1)^m t^{2m}}{(2m)!} = \cos^{-n}(t) = \sec^n(t). 
\]
This is the exponential generating function of the sequence $T_{m,i}$~\cite[A085734]{OEIS}.
Therefore we have the desired expression of $E_{2m}^c(n)$ as a polynomial in $n$.

On the other hand, the odd powers $t$ in $e^t \cosh^{n-1}(t)$ are given by $\sinh(t)\cosh^{n-1}(t)$, and the derivative of $\cosh^n(t)$ is $n\sinh(t)\cosh^{n-1}(t)$.
Thus $E_{2m}^c(n) = n E_{2m-1}^c(n)$.
\end{proof}

\begin{remark}\label{rem:T}
By Corollary~\ref{cor:T}, $E_{2m}^c(n)$ is a polynomial in $n$ with degree $m$, leading coefficient $(2m-1)!!$, and constant term given by the \emph{tangent number}~\cite[A000182]{OEIS}. 
It agrees with the zeta polynomial for the poset of even-sized subsets of $\{1,2,\ldots,2m\}$ ordered by inclusion by a comment of Geoffrey Critzer (Apr 22 2023) in OEIS~\cite[A085734]{OEIS}.
Another comment by David Callan (Oct 25 2004) in OEIS~\cite[A085734]{OEIS} asserts that $T_{m,k}$ enumerates down-up permutations on $\{1,2,\ldots,2m+2\}$ with $k+1$ left-to-right maxima. 
For example, $T_{1,1}$ counts the following $3$ down-up permutations on $\{1,2,3,4\}$ each with $2$ left-to-right maxima: $2143$, $3142$, $3241$.
Also, $T_{m,k}/2^{m-k}$ is an integer with various properties~\cite[A083061, A094665]{OEIS}.
It would be interesting to explore connections between these results and our definition of $E_{2m}^c(n)$.
\end{remark}

\section{Proof of Theorem~\ref{thm4}: combinatorial interpretation of $E^c_m(n)$}\label{sec:interpretation}

In this section we build a coin-flipping model to interpret $E_m^c(n)$ and the more general number $c_{m,k}$ as defined in Lemma~\ref{lem:c}.

\begin{proposition}
If we begin with $n-1$ coins all showing tails and flip at most one of them at each step, then the number of ways to have exactly $k$ tails right after $m$ steps is the number $c_{m,k}$ defined in Lemma~\ref{lem:c}.
\end{proposition}

\begin{proof}
Let $c'_{m,k}$ be the number of ways to have exactly $k$ heads right after $m$ steps.
One can check that $c'_{0,0}=1$ and $c'_{m,k}=0$ whenever $k<0$ or $k>m$.
For $m\ge1$, we have
\[
c'_{m,k} = (n-k) c'_{m-1,k-1} + c'_{m-1,k} + (k+1)c'_{m-1,k+1}
\]
since there are exactly $k$ heads at step $m$ if and only if one of the following occurs.
\begin{itemize}
\item
There are $k-1$ heads at step $m-1$, and one of the $n-k$ tails is flipped at step $m$.
\item
There are $k$ heads at step $m-1$, and nothing is changed at step $m$.
\item
Then are $k+1$ heads st step $m-1$, one of which is flipped at step $m$.
\end{itemize}
It follows that $c'_{m,k} = c_{m,k}$.
\end{proof}

Taking $k=0$ gives Theorem~\ref{thm4}: beginning with $n-1$ coins all showing tails and flipping at most one of them at each step, the number of ways to have all tails right after $m$ steps is $c_{m,0} = E^c_m(n)$.

\begin{remark}
It seems interesting to study $c_{m,k}$ as a polynomial in $n$ for $k\ge1$, although we cannot find its coefficients in OEIS~\cite{OEIS} (see Table~\ref{tab4} for small values of $m$ and $k$).
Note that $c_{m,k} = (n-1)_k d_{m,k}$ by Remark~\ref{rem:d}.
\begin{table}[h]
\tiny
\begin{tabular}{c|ccccccccc}
\diagbox{$m$}{$k$}& $0$ & $1$ & $2$ & $3$  \\
\hline
$0$ & $1$ & \\ 
$1$ &  $1$ & $n-1$ \\ 
$2$ &  $n$ & $2n-2$ & $n^2-3n+2$ \\ 
$3$ &  $3n-2$ & $3n^2-5n+2$ & $3n^2-9n+6$ & $n^3-6n^2+11n-6$ \\ 
$4$ &  $3n^2-2n$ & $12n^2-28n+16$ & $6n^3-26n^2+36n-16$ & $4n^3-24n^2+44n-24$  \\ 
$5$ &  $15n^2-30n+16$ & $15n^3-45n^2+46n-16$ & $30n^3-150n^2+240n-120$ & $10n^4-80n^3+230n^2-280n+120$ \\ 
$6$ & $15n^3-30n^2+16n$ & $90n^3-390n^2+572n-272$ & $45n^4-285n^3+676n^2-708n+272$ & $60n^4-520n^3+1620n^2-2120n+960$ 
\end{tabular}
\vskip5pt
\caption{$c_{m,k}$ for small values of $m$ and $k$}
\label{tab4}
\end{table}
\end{remark}

\section*{Acknowledgment}

The authors used SageMath to help discover and verify some of the results in this paper.


\end{document}